\documentclass[twoside,10pt]{article}
\usepackage{amssymb}

\newtheorem{thm}{Theorem}

\def \k {\kern 5 pt}
\def \n {\noindent}

\usepackage{amssymb}

\def \k {\kern 5 pt}
\def \n {\noindent}

\title{\bf Quasi-invariant measures on the path space of a diffusion}
\date{}
\begin{document}
\maketitle

\kern - 30 pt

\centerline { Denis Bell\footnote{Research partially 
supported by
NSF grant DMS-0451194.}}

\centerline {Department of Mathematics, University
 of North Florida,}

\centerline {4567 St. Johns Bluff Road South, Jacksonville, FL 32224,
U. S. A.}
\centerline {email: dbell@unf.edu, phone: 904-620-2653, fax: 904-620-2818.} 

\kern 50 pt

{\bf Abstract} 

The author has previously constructed a class of admissible vector fields
 on the path space of an elliptic diffusion process $x$ taking values	 in a closed compact
 manifold. 
In this Note the existence of flows for this class of  vector fields is established and it is 
shown that
 the law of $x$ is quasi-invariant under these flows.

\kern 20 pt

\n {\bf  R\'esum\'e}

L'auteur a pr\'ec\'edemment construit une classe de champs de
vecteurs admissibles 
 sur l'espace des chemins d'une diffusion elliptique $x$ prenant valeurs dans une 
vari\'et\'e compacte ferm\'ee. 
Dans cette Note l'existence des flots pour cette classe de champs de vecteurs 
est \'etablie et on   montre que
 la loi de $x$ est quasi-invariante sous ces flots.

\vfil\break

 {\bf 1. Introduction}

Let $N$ denote a smooth manifold  equipped
with a finite Borel measure $\nu$. A vector field $Z$ on $N$ is said to be
{\it admissible} (with respect to $\nu$) if there exists an $L^1(\nu)$ random variable
 $Div(Z)$ such that the equality
$$E\big [Z(\Phi)(x)\big ] = E\big [\Phi(x)Div(Z)\big ]$$
holds for a dense class of real-valued functions $\Phi$ on $N$. Suppose $\{x\mapsto x^s, s \in {\bf R}\}$ is the flow
on $N$ generated by  $Z$. Assume  $x^0$ is a random variable with law $\nu$ and
let $\nu^s$ denote the law of $x^s$. Then $\nu$ is said to be {\it
quasi-invariant} under  the flow of $Z$ if $\nu_s$  and $\nu$ are 
equivalent measures, 
for all $s$. 

These two properties are closely related. For example, suppose $\nu$ is quasi-invariant under
the flow of a vector field $Z$  
and write $\rho_s$ for the family of Radon-Nikodym derivatives 
$
\rho_s={d\nu_s\over d\nu}.
$
Then an obvious calculation suggests that
 $Z$ will generally be
admissible, with 
$$
Div(Z) = {d\rho_s\over ds}/_{s=0}.
$$
The converse is not necessarily true, as is easily seen by considering
the case where $\nu$ is a measure on a Euclidean space with a smooth compactly
supported density and $Z$ is a constant vector field. However, admissibility of $Z$ with
 respect to $\nu$, together with the existence
 of a suitably regular version of the process $Div(Z)(x^s)$,  can be shown to imply 
quasi-invariance
of $\nu$ under the flow $x^s$ of $Z$ (cf. Theorem 3 below).  

In this Note we study these properties  in the  following
 setting. 
Let  $X_1,\dots, X_n$ and $Y$ denote smooth vector fields defined on a
closed compact $d$-dimensional manifold $M$. Consider the Stratonovich
stochastic differential equation (SDE) 
$$
dx_t=\sum_{i=1}^nX_i(x_t)\circ dw_i + Y(x_t)dt, \k t \in [0,T] \eqno(1)
$$
$$
x_0=o \kern 155 pt
$$
where $(w_1,\dots,w_n)$ is a  standard Euclidean Wiener process and $o$ is a
point in $M$. Assume the diffusion (1) is {\it elliptic}, i.e. the
vector fields $X_1,\dots, X_n$ span $TM$ at each point of $M$. 
Define $\nu$ to be the law of $x$, considered as a
measure on  $C_o(M)$, the space of continuous paths 
$\{\sigma: [0,T]\mapsto M/  \kern 2 pt
\sigma (0)=o\}$.  

The author  ([1] and [2]) has recently constructed a new class of 
admissible  vector fields on 
the path space $\big (C_o(M), \nu\big )$. These vector fields, and 
their divergences, are described in Theorem 1.
Theorem 2 asserts  that the vector fields in Theorem 1
generate  flows on $C_o(M)$. Theorem 3 gives a criterion under which admissibility 
of a vector field
 implies quasi-invariance of the corresponding flow of measures.
 The main result in the paper is Theorem 4, which 
asserts that  the measure $\nu$  is quasi-invariant under the  
flows generated by the vector fields in  Theorem 1.

The issue of quasi-invariance of measures under flows of vector fields  
 has been well-studied
 in the
the past two decades. Significant results in this area have been obtained by Driver
 [5] and Hsu [6] in
the classical path space framework studied here and by Cruzeiro [4] for vector
fields on abstract Wiener space. The results
presented here are a continuation  of this tradition.

\k

\k

\k

{\bf 2. Additional notations}

Denote by $[g_{jk}]_{j,k=1}^d$ the
Riemannian metric on $M$ defined by $g^{jk} = a_{ij}a_{ik}$ where $X_i= a_{ir}\partial/\partial 
x_r$ is the local representation of $X_i, 1\le i\le n$ (here and henceforth, we adopt the
usual summation
convention: whenever an index in a product is repeated, that index is assumed to be summed on). Let
 $\nabla$
denote the Levi-Civita covariant derivative associated with this metric.

Define a set of 1-forms $\omega^{jk}, 1\le j, k \le n$ on $M$ by
$$\omega^{jk}(.)= <\nabla_{X_j}X_k,.>-<\nabla_.X_j,X_k> 
$$
and functions $B^{jk}, 1\le j, k \le n$ on $M$ by
$$
B^{jk}(x)= {1\over 2}\Big (<L_{ji}X_i, X_k>-<L_{ij}X_k, X_i>  - <\nabla_{X_j}X_k,
 \nabla_{X_i}X_i> \kern 20 pt$$
$$\kern 160 pt  +<\nabla_{X_p}X_i,X_k><\nabla_{X_j}X_p, X_i>
 \Big) (x)
$$
where 
$L_{ij}$ is the differential operator $\nabla_{X_i}\nabla_{X_j}-\nabla_{\nabla_{X_i}{X_j}}. 
$

\k

\k

\k

{\bf 3. Statement of results}

The methods of [1] and [2] yield the following 
\begin{thm} 
Let $r=(r^1,\dots, r^n)$ be any path in the Cameron-Martin space of ${\bf R}^n$ and 
define
 $h^i, 1\le i\le n$ by the
following system of SDE's
$$
dh^i_t=  \omega^{ji}(\circ dx_t)h^j_t+\big [B^{ji} +
 <\nabla _{X_j}Y,X_i>\big ](x_t)h^j_tdt + \dot r^i_tdt
$$
$$
h^i_0=0. \kern 220 pt
$$
Then the vector field $Z(x)_t \equiv X_i(x_t)h^i_t, \kern 2 pt t \in 
[0, T]$ on $C_o(M)$ is admissible and
$$Div(Z)=\int_0^T\Big (\dot r^i_t+{1\over 2}<Ric(Z_t),X_i(x_t)>\Big)dw_i
$$
where $Ric$ denotes Ricci curvature.
\end{thm}
Although the equations defining $B^{jk}$ are lengthy, the {\it form} of the vector 
field $Z$ is 
relatively simple. The next result establishes the existence of flows on $C_o(M)$
for vector fields of this form.

\begin{thm} Let $T^{ij}$ and $f^{ij}, 1\le i, j\le n$ be, respectively, 
smooth 1-forms ond real-valued functions on $M$ and suppose $g^i \in L^2[0,T], 
1\le i\le n$ are deterministic real-valued functions. Define a vector field
$V$ on $C_0(M)$ by $V_t=X_i(x_t)\eta^i_t$, where $\eta_1,
\dots,\eta_n$ solve the  
 system of SDE's
$$
d\eta^i_t =T^{ij}(\circ dx_t)\eta^j_t+\big [f^{ij}(x_t)\eta^j_t + g^i(t)\big ]dt,
\k t \in [0,T]$$
$$
\eta^i_0=0. \kern 180 pt
$$

\goodbreak
Then there exists a solution $x^s$ in $C_o(M)$ to the flow equation
$$
{dx^s\over ds} = V(x^s), \kern 5 pt s\in {\bf R} \eqno (2)
$$
$$
x^0= x  \kern 55 pt
$$
where $x$ is the process in equation (1).

Furthermore, the paths $x^s$ are semimartingales of the  form
$$
dx^s_t =X_i(s, t)\circ dw_i + X_0(s, t)dt \eqno (3)
$$
where $X_j(s,.), 0\le j\le n$ are adapted processes (with respect to $w$)
 in $TM$ such that, on a set of 
full Wiener 
measure, 
$s\mapsto X_j(s,.)$ is continuous  into the space $L^2{[0,T]}$.
\end{thm}

We prove Theorem 2  using a version of the Picard iteration scheme. 
The proof draws upon the closure
of the class of It\^o processes under composition with smooth maps, together
with standard  estimates on stochastic integrals. Each $m$-th stage of the iteration 
 yields a family of It\^o
 processes $dx^{s,m}_t =X_i^m(s, t)\circ dw_i + X_0^m(s, t)dt$. Expressions 
are obtained 
for the coefficients $X_j^m, 0\le j\le n$, in terms of $\{X_j^{m-1}, 0\le j\le n\}$.
These expressions are used to prove  that
$X_j^m(s,t)$ converge to limits $X_j(s,t)$ that have the stability in $s$ 
described in
the Theorem.  The iterative procedure is then shown to imply that 
the processes $x^s$ defined by (3) satisfy the flow equation (2).

\k

\n The following result is proved in   [3]

\begin{thm} Suppose $N$ is a manifold equipped with a finite Borel measure
 $\gamma$. Let
$Z$ be an admissible vector field on $N$ with flow $\{ \sigma_s, s
 \in {\bf R}\}$ and let $\gamma_s$
 denote the measure $\sigma_s(\gamma)$. Suppose there exists 
$B\subseteq N$ with $\gamma_s(B) = 1$ for all $s$, such that Div(Z) is defined
 and $Z$-differentiable on $B$, and that the function $s\mapsto 
Div(Z)\big (\sigma_s(p)\big )$ is absolutely continuous for $p\in B$. Then 
$\gamma_s$ and $\gamma$ are equivalent measures and
$$
{d\gamma_s\over d\gamma}(p) = \int_0^sDiv (Z)\big (\sigma_{-u}(p)\big )du
$$ 

\end{thm} 
Combining Theorems 1, 2, and 3, we obtain our main result

\begin{thm} Let $x^s$  denote the flow  
in Theorem 2 generated by the  vector field $Z$ in Theorem 1. Then for each $s$, the law $\nu_s$ of
$x^s$ is equivalent to  $\nu$ and
$${d\nu_s\over d\nu}(x)= \exp\int_0^sDiv (Z)(x^{-u})du.
$$
\end{thm}

\k

{\it Remark}. Theorem  4 is similar to a result in Bruce Driver's
paper [5]. However, the methods are different. Driver studied the class of vector fields 
on $C_o(M)$
 obtained from Cameron-Martin paths 
in 
$T_oM$ by stochastic parallel translation along the paths of the diffusion process $x$. He 
constructed   
quasi-invariant flows for these vector fields, then
 used the quasi-invariance property 
 to
deduce admissibility of the vector fields. In the work described here, this approach 
is reversed. We first construct a class of admissible vector fields, then show
that they generate flows and quasi-invariant measures on $C_o(M)$. 

In [7], Hu, \"Ust\"unel and Zakai   
studied rotation type measure-preserving transformations of abstract Wiener space 
and constructed flows associated to these transformations. In contrast to the flows studied here and in
Driver's work, the transformations considered 
in [7] are non-adapted (the notion of adaptedness having no obvious meaning
 in the context of an 
abstract Wiener space).

\k

\k

\k

\centerline {\bf References}

\k

 \n [1] D. Bell,      Divergence theorems in path space. {\it J. Funct. Anal}.
{\bf 218} (2005), no. 1, 130 - 14.

\n [2] D. Bell, Divergence theorems in path space II: degenerate diffusions, 
{\it C. R. Acad. Sci. Paris Sïr. I Math.}, to appear.

\n [3] D. Bell, Admissible vector fields and quasi-invariant measures.
 Appendix to
 {\it The Malliavin Calculus}, 2nd edition. Dover Publications, Mineola, 
NY, 2006.

\n [4] A. B. Cruzeiro,  \' Equations diff\' erentielles sur l'espace de Wiener et formules de
 Cameron-Martin non-lin\' eaires. {\it
J. Funct. Anal.} {\bf 54} (1983) 206-227.

\n  [5] B. Driver, A Cameron-Martin type quasi-invariance
theorem for Brownian motion on a compact manifold. {\it
J. Funct. Anal.} {\bf 109} (1992) 272-376.

\n [6] E. P. Hsu,  Quasi-invariance of the Wiener measure on the path 
space over a 
compact
 Riemannian manifold.  {\it
J. Funct. Anal.} {\bf 134} (1995) 417-450.

\n [7] Y. Hu, A. S. \"Ust\"unel and M. Zakai, Tangent processes on Wiener 
space. 
{\it J. Funct. Anal.} {\bf 192} (2002), no. 1, 234-270.

\goodbreak

\end{document}